\documentclass{amsart}
\usepackage[subpreambles=true]{standalone}
\usepackage{amsmath,amssymb,nicefrac,import,placeins}

\newcommand{\romannum}[1]{\uppercase\expandafter{\romannumeral #1\relax}}

\begin{document}

\title{The combinatorial method to compute the sum of the powers of primes}
\author[A. Orlov] {Alexey Orlov}
\subjclass[2010]{11N05,11Y16,11Y55,11Y70}

\begin{abstract}
We will generalize the combinatorial algorithms for computing $\pi(x)$ to compute sums ${F(x) = \sum_{p \leq x} p^k}$
for $k \in \mathbb{Z}_{\geq 0}$.
The detailed exposition of algorithms is included along with implementation details.
\end{abstract}
\maketitle
\tableofcontents

\section{Introduction}

The history of the calculation of $\pi(x)$ without enumerating all the primes up to $x$ dates back to Legendre. His
method was combinatorial at heart and was based on the Inclusion-Exclusion Principle, requiring primes up to $\sqrt{x}$.
The next step was done by Meissel in 1870, and later by Lehmer \cite{Lehmer} in 1959,
who streamlined the Meissel's algorithm into what would be later called Meissel-Lehmer algorithm.
Next improvements came from Lagarias, Miller, Odlyzko \cite{LagariasMillerOdlyzko} in 1985, and from
Deleglise, Rivat \cite{DelegliseRivat} in 1996. Later improvements were mainly concerned with implementation
improvements and are of no interest to us here. All these algorithm are ``combinatorial'', as opposed to ``analytical''
methods that were inspired by the works of Lagarias, Odlyzko \cite{LagariasOdlyzkoAnal}.
Our goal is to generalize the ``combinatorial'' method to calculate not just $\pi(x) = \sum_{p \leq x}p^0$, but the
more general sum $F_k(x) = \sum_{p \leq x}p^k$ for $k$ non-negative integer. \\
We will give a self-contained exposition of the ``combinatorial'' algorithm for calculating $F_k(x)$ and provide
some values of $F_k(x)$ for $k=2,3,4$.

\ifstandalone
\bibliography{master}
\bibliographystyle{amsplain}
\fi

\section{General Description}

We shall start with describing the combinatorial method, closely following \cite{Lehmer}.
Let
\[F(x) = \sum_{p \leq x}f(p),\]
where $f$ is completely multiplicative. We write
\[m_a=\prod_{i \leq a}{p_i}\]
for the product of the first $a$ primes. \\
By $\phi(x,a)$ we denote the sum of $f(n)$ over the numbers $n \leq x$ divisible by none of the first $a$ primes:
\[\phi(x,a) = \sum_{\substack{n \leq x \\ (n,m_a) = 1}}{f(n)};\]
and by $P_k(x,a)$ we denote the sum of $f(n)$ over $n \leq x$ such that $n$ is a product of $k$ primes
each greater than $p_a$:
\[P_k(x,a) = \sum_{p_a < q_1, \dots, q_k \leq x}{f(q_1 \dots q_k)}.\]
Also, by the usual convention, we set $P_0(x, a) = f(1) = 1$. \\
We clearly see that
\[\phi(x,a) = \sum_{k=1}^\infty P_k(x,a).\]
Further, for $k$ such that $x < p_{a+1}^k$ we have $P_k(x,a) = 0$, and upon writing $r$ for the smallest such $k$,
we may rewrite this sum as
$\phi(x,a) = \sum_{k=0}^{r-1}{P_k(x,a)}$. \\
We can also go in the ``opposite'' direction: fixing $r$, \ie limiting the number of $P_k$, we will acquire
the bounds on $a$:
\[a \in \left[\pi\left(x^{\frac{1}{r}}\right), \pi\left(x^{\frac{1}{r-1}}\right)\right).\]
We fix a parameter $Y$ such as
\[x^{\frac{1}{r}} \leq Y < x^{\frac{1}{r-1}},\]
and set $a = \pi(Y)$. We will return to the choice of $Y$ later, when we will use it to balance different parts of
computation to achieve optimal performance. \\
After we expand $P_1$:
\[P_1(x, a) = \sum_{p_a < p \leq x}{f(p)} = F(x) - F(p_a),\]
we may rewrite the sum as
\[F(x) = \phi(x,a) + F(p_a) - 1 - \sum_{k=2}^{r-1}{P_k(x,a)},\]
thus the computation of $F$ is reduced to the computation of $P_k$ and $\phi$. \\

We shall start with the computation of $\phi$. We define $Q(x,k) = \sum_{ik \leq x} {f(ik)}$
and we shall prove that $\phi(x,a) = \sum_{d \mid m_a}{\mu(d)Q(x,d)}$.
\begin{align*}\sum_{d \mid m_a}{\mu(d)Q(x, d)}
    &= \sum_{d \mid m_a}{\mu(d)} \sum_{\substack{d \mid m \\ m \leq x}}{f(m)}
     = \sum_{m \leq x}{f(m)} \sum_{d \mid (m, m_a)}{\mu(d)}                     \\
    &= \sum_{\substack{m \leq x \\ (m,m_a) = 1}}{f(m)} = \phi(x,a).
\end{align*}
We now split all the divisors $d$ of $m_a$ in two groups. \\
The contribution of $d$ such that $p_a \nmid d$ is
\[\sum_{\substack{d \mid m_a, p_a \nmid d}}{\mu(d)Q(x, d)} = \sum_{d \mid m_{a-1}}{\mu(d)Q(x, d)} = \phi(x, a-1),\]
and for the rest
\begin{align*}\sum_{\substack{p_a \mid d \\ d \mid m_a}}{\mu(d)Q(x, d)}
    &= \sum_{d \mid m_{a-1}}{\mu(p_a d)Q(x, p_a d)}
     = -\sum_{d \mid m_{a-1}}{\mu(d)} \sum_{\substack{p_a d \mid m, m \leq x}}{f(m)}                        \\
    &= -f(p_a) \sum_{d \mid m_{a-1}}{\mu(d)} \sum_{\substack{d \mid m, m \leq \nicefrac{x}{p_a}}}{f(m)}     \\
    &= -f(p_a) \phi(\nicefrac{x}{p_a}, a-1),                                                                \\
\end{align*}
giving us
\[\phi(x,a) = \phi(x,a-1) - f(p_a)\phi\left(\frac{x}{p_a}, a-1\right).\]
As $\phi(x,0) = \sum_{n \leq x}f(n)$ we see that this recursive formulation may be represented as a binary tree
of the height $a$.

Next we shall expand $P_2$ and $P_3$:
\begin{align*} P_2(x,a)
    &= \sum_{\substack{p_a < p_i,p_j \\ p_i p_j \leq x}}f(p_i p_j)
        = \sum_{a < i \leq \pi(\sqrt{x})}f(p_i)\sum_{i \leq j \leq \pi(\nicefrac{x}{p_i})}f(p_j)    \\
    &= \sum_{a < i \leq \pi(\sqrt{x})}f(p_i)\left[F\left(\frac{x}{p_i}\right) - F(p_i-1)\right].    \\
P_3(x,a)
    &= \sum_{\substack{p_a < p_i,p_j,p_k \\ p_i p_j p_k \leq x}}f(p_i p_j p_k)                      \\
    &= \sum_{a < i \leq \pi(\sqrt[3]{x})}f(p_i)
            \sum_{i \leq j \leq \pi\left(\sqrt{\nicefrac{x}{p_i}}\right)}f(p_j)
            \sum_{j \leq k \leq \pi\left(\nicefrac{x}{p_i p_j}\right)}f(p_k)                        \\
     &= \sum_{a < i \leq \pi(\sqrt[3]{x})}f(p_i)
            \sum_{i \leq j \leq \pi\left(\sqrt{\nicefrac{x}{p_i}}\right)}
            f(p_j)\left[F\left(\frac{x}{p_i p_j}\right) - F(p_j-1)\right].
\end{align*}
We see that the sum for $P_3$ is way more complicated than $P_2$, this is the price we pay to lower $a$ from
$\sqrt[3]{x}$ to $\sqrt[4]{x}$. Thankfully, as we will see later, in the case of $r=3$ we can optimize the computation
of $\phi(x,a)$ in such a way that we will need the ``recursive'' procedure above only for $\phi(x,\sqrt[4]{x})$, thus
making the addition of $P_3$ unneeded. \\
Hence, from now on, we fix $r=3$, and we consider $Y$ such that $\sqrt[3]{x} \leq Y < \sqrt{x}$, to be a parameter,
setting $a = \pi(Y)$. We shall now consider the computation of $\phi(x,a)$ and $P_2$ in greater detail.

\ifstandalone
\bibliography{master}
\bibliographystyle{amsplain}
\fi

\section{Computing $P_2(x,a)$}

We recall the formula for calculating $P_2$:
\begin{align*}
P_2(x,a)
        &= \sum_{p_a < p \leq \sqrt{x}}\left[F\left(\frac{x}{p_i}\right) - F(p_i-1)\right]     \\
        &= \sum_{p_a < p \leq \sqrt{x}}f(p)F\left(\frac{x}{p_i}\right) - \sum_{p_a < p \leq \sqrt{x}}f(p)F(p_i-1).
\end{align*}
We write $P_2(x,a) = S_1 - S_2$, where
\begin{align*}
S_1 &= \sum_{a < i \leq \pi(\sqrt{x})}f(p_i)F\left(\frac{x}{p_i}\right) \\
S_2 &= \sum_{a < i \leq \pi(\sqrt{x})}f(p_i)F(p_i-1)
\end{align*}
For calculating $S_2$, we note that we can accumulate $F(p_i-1)$ as we go through primes in
$\left(p_a,\sqrt{x}\right]$. \\
To calculate $S_1$ we need to know the values of $F\left(\nicefrac{x}{p_i}\right)$ and we have
\[\sqrt{x} \leq \frac{x}{p_i} \leq \frac{x}{p_{a+1}} < \frac{x}{Y}.\]
We write
\[L = \left\lceil \sqrt{x} \right\rceil, R = \left\lceil \frac{x}{Y} \right\rceil\]
and proceed by sieving the interval $\left[L, R\right)$ in blocks $I_k = \left[L+(k-1)B, L+kB\right)$ of size $B$,
with the possible exception of the last block which can be shorter. As we have $\frac{x}{p_i} \in I_k$, it is natural
to consider the block
\[J_k = \left(\frac{x}{L+kB}, \frac{x}{L+(k-1)B}\right] \cap \left(Y,\sqrt{x}\right],\]
and sieve it completely.
In the end we know all the primes $p_i$ such that $\frac{x}{p_i} \in I_k$. As we have sieved $I_k$ fully we can
calculate $F\left(\frac{x}{p_i}\right)$ easily,
so we just add $f(p_i)F\left(\nicefrac{x}{p_i}\right)$ to the $S_1$ accumulator. \\
We should estimate the maximal length of $J_k$. Suppose we consider $I_k = [h-B,h)$, with $h > \sqrt{x}$. Then
\[\frac{x}{h-B} - \frac{x}{h} = \frac{xB}{h(h-B)}.\]
If we assume that $h-B \geq \sqrt{x}$ (that is the usual block for $S_1$ calculation) then $h(h-B) \geq x$,
and we clearly don't need more than $B$ integers. \\
On the other hand, if $h-B < \sqrt{x}$, we would have $J_k$ shortened to
$\left(\frac{x}{h}, \sqrt{x}\right]$ and then
\[\sqrt{x}-\frac{x}{h} = \frac{h\sqrt{x}-x}{h} \leq \frac{h\sqrt{x}-x}{\sqrt{x}} = h - \sqrt{x} < B,\]
and again, we don't need more than $B$ integers. \\

\ifstandalone
\bibliography{master}
\bibliographystyle{amsplain}
\fi

\section{Computing $\phi(x,a)$}

\subsection{Recursion Tree}

It is clear that as we go through the recursion
\[\phi(x,a) = \phi(x,a-1) - f(p_a)\phi(\nicefrac{x}{p_a}, a-1),\]
we obtain a binary tree, and ultimately we need to sum the leaf values.
The nodes of the tree are of the form
\[\mu(n)f(n)\phi\left(\frac{x}{n},b\right) \text{, where } n = p_{i_1} \dots p_{i_r}, a \geq a_1 > \dots > a_r > b,\]
and we will label them with $(n,b)$. \\

Now, we need to consider the possibility of the ``truncation'' of our tree. We will consider the case when $b$ is
``small'', and when $n$ is ``large'' (\ie $\nicefrac{x}{n}$ is ``small'').

The very first truncation rule is quite obvious: for $z \leq p_k$ we have $\phi(z,k) = 1$, and thus we have \\
\textbf{The Truncation Rule $T_0$} Stop at the node $(n,b)$ if either of the following holds:
\begin{enumerate}
  \item $n \geq \nicefrac{x}{p_b}$
  \item $n < \nicefrac{x}{p_b}$, and $b = 0$.
\end{enumerate}

We shall write $\phi_m$ for the $\phi$ function corresponding to $f(n) = n^m$. In \cite{Lehmer} it was proposed to
tabulate the values of $\phi_0(x,a)$ for small values of $a$ when computing $\pi(x)$, and we shall generalize this.
Let $P$ be the product of the first $K$ primes. Now we shall consider the task of calculating $\phi_m(x,K)$, having
\[x = qP + r \text{ with } 0 \leq r < P.\]
Assuming that we have precomputed $\phi_m(n,K)$ for $n \leq P$ we use the Binomial theorem to obtain
\begin{align*} \phi_m(x,K)
    &= \sum_{i=0}^{q-1}\sum_{\substack{j < P \\ (j,P) = 1}}(iP + j)^m
        + \sum_{\substack{j \leq r \\ (j,P) = 1}}(qP + j)^m                                 \\
    &= \sum_{i=0}^{q-1}\sum_{\substack{j < P \\ (j,P) = 1}}\sum_{k=0}^m\binom{m}{k}i^kP^kj^{m-k}
        + \sum_{\substack{j \leq r \\ (j,P) = 1}}\sum_{k=0}^m\binom{m}{k}q^kP^kj^{m-k}      \\
    &= \sum_{k=0}^m\binom{m}{k}P^k \sum_{i=0}^{q-1}i^k \sum_{\substack{j < P \\ (j,P) = 1}}j^{m-k}
        + \sum_{k=0}^m\binom{m}{k}q^kP^k\sum_{\substack{j \leq r \\ (j,P) = 1}}j^{m-k}      \\
    &= \sum_{k=0}^m\binom{m}{k}P^k \phi_{m-k}(P-1,K) \sum_{i=0}^{q-1}i^k
        + \sum_{k=0}^m\binom{m}{k}q^kP^k \phi_{m-k}(r,K)                                    \\
    &= \sum_{k=0}^m\binom{m}{k}P^k
        \left(\phi_{m-k}(P-1,K) \sum_{i=0}^{q-1}i^k + \phi_{m-k}(r,K)q^k\right).
\end{align*}
We may assume, that we can calculate $\sum_{i=0}^{q-1}i^k$ for a fixed $k$ efficiently, and then we need to store
the values of $\phi_0(x,K), \dots \phi_m(x,K)$ for $x < P$. Thus we can update our truncation rule: \\
\textbf{The Truncation Rule $T_1(K)$} Stop at the node $(n,b)$ if either of the following holds:
\begin{enumerate}
  \item $n \geq \nicefrac{x}{p_b}$
  \item $n < \nicefrac{x}{p_b}$, and $b = K$.
\end{enumerate}
We shall briefly note that setting $K=0$ gives us the initial truncation rule. \\
Now we will update the truncation rule for ``large'' $n$, following \cite{LagariasMillerOdlyzko}: \\
\textbf{The Truncation Rule $T_2(Y,K)$} Stop at the node $(n,b)$ if either of the following holds:
\begin{enumerate}
  \item $n > Y$
  \item $n \leq Y$, and $b = K$
\end{enumerate}
We will call the leaves of type 1 \textit{special leaves} and of type 2 \textit{ordinary leaves}. \\
At this point we should prove the correctness of this rule, \ie that we account for all the nodes, and exactly once.
Consider the level of the tree corresponding to $K$, consisting of the nodes $(n,K)$.
Since we go through the primes in descending order we clearly have $(P,n) = 1$,
namely $n$ is not divisible by the first $K$ primes. On this level we have ordinary nodes $(n,K)$ with
$n \leq Y$, and the nodes with $n > Y$, which can be backtracked to the special node. \\
Therefore we have
\[
\phi(x,a) = \sum_{(n,K) \text{ ordinary}} \mu(n)f(n)\phi\left(\frac{x}{n},K\right)
    + \sum_{(n,b) \text{ special}} \mu(n)f(n)\phi\left(\frac{x}{n},b\right).
\]
We note that the contribution of the ordinary leaves can be computed immediately, as it is the sum over squarefree
$n \leq Y$ such that $(n,P)=1$, and the terms can be computed efficiently.
Thus the special leaves make up the core essence of the calculation. \\

\subsection{Special Leaves}

For the special leaf $(n,b)$ we note that its parent couldn't be $(n,b+1)$ as we would stop earlier. Hence it was
$(n^\ast, b+1)$, with $n = n^\ast p_{b+1}$, and we have $n^\ast \leq Y < n^\ast p_{b+1}$, $l_p(n) = p_{b+1}$,
where we write $l_p(n)$ for the smallest prime factor of n; as we go trough the primes in descending order we must
have $l_p(n^\ast) > p_{b+1}$. Thus the multiplier for the node $(n,b)$ is
\[
    \mu(n)f(n) = \mu(n^\ast)\mu(p_{b+1})f(n^\ast)f(p_{b+1}) = -\mu(n^\ast)f(n^\ast)f(p_{b+1}).
\]
We might go through the primes $p_k$ with $k \in [1,a]$ and enumerate squarefree
$n^\ast \in \left(\nicefrac{Y}{p_k}, Y\right]$, with $l_p(n^\ast) > p_k$, calculating $\phi$ recursively. \\

We shall now describe the procedure to compute the contribution of special leaves without recursion.
We note that we can actually compute $\phi(x,a)$ iteratively in blocks. Consider the block $I = [l,h)$ of length
$B = h-l$ and suppose that we know $\phi(l-1,k)$ for $0 \leq k \leq a$. For $x \in I$ and $k \leq a$ we have
\[\phi(x,k) = \phi(l-1,k) + \sum_{\substack{n \in I, n \leq x \\ (n,m_k) = 1}}f(n).\]
Say we have a zero-based array
\[f_{block} = [f(l+0), f(l+1), \dots f(l+B-1)].\]
Then for $k$ from $0$ to $a$ we know that
\[\phi(x,k) = \phi(l-1,k) + \sum_{i=0}^{x-l}f_{block}[i],\]
and after we have calculated all the $\phi$ values we need, we ``strike out'' the multiples of $p_k$, \ie setting
$f_{block}[i]$ to zero for $p_k \mid l+i$, and proceed to the next $k$. \\
Thankfully, there is an efficient data structure for calculating prefix sums of the mutable array: the Fenwick tree.
Thus we may assume that this can be done efficiently. \\

For a special leaf $(n,b)$ we have $n > Y$ and $\nicefrac{x}{n} < \nicefrac{x}{Y}$. We will sieve interval
$\left[1, \nicefrac{x}{Y}\right)$ in blocks of size $B$
\[I_k = \left[1+(k-1)B, 1+kB\right) = \left[1+(k-1)B, kB\right],\]
with the last block potentially shortened. \\
Since $n = n^\ast p_{b+1}$, we can process this leaf after we have sieved the first $b$ primes.
Thus we get bounds
\[\frac{x}{((k+1)B+1)p_{b+1}} < n^\ast \leq \frac{x}{(kB+1)p_{b+1}},\]
and
\[
    n^\ast \in \left(\frac{x}{((k+1)B+1)p_{b+1}}, \frac{x}{(kB+1)p_{b+1}}\right] \cap \left[1, Y\right],
\]
such that $l_p(n^{\ast}) > p_{b+1}$, and $\mu(n^{\ast}) \not = 0$. \\
We recall that a special node $(n,b)$ is such that $(n,P) = 1$, thus for the first $K$ primes we only need to
``strike out'' their multiples; all this machinery is needed only starting with $b = K+1$. We can even reuse
sieving results obtained from the precalculation of $\phi$, if this proves efficient. \\

Now, we proceed to lower the number of primes considered per block, from the first $a$ primes to the first
$\pi\left(\sqrt[4]{x}\right)$ primes. Following \cite{DelegliseRivat} we should further split special leaves.
We will consider the 3 leaf classes:
\begin{enumerate}
  \item $\sqrt[3]{x} < l_p(n) \leq Y$
  \item $\sqrt[4]{x} < l_p(n) \leq \sqrt[3]{x}$
  \item $l_p(n) \leq \sqrt[4]{x}$
\end{enumerate}
We also note that if we use the precomputed $\phi$ table, the special leaves $(n,b)$ must have $b > K$. We might tweak
lower bounds in the above, but this is way too complicated. It is way easier to either lower $K$ to have
$K < \pi\left(\sqrt[4]{x}\right)$ or to compute $F(x)$ directly (as $x$ is quite low in this case).

First of all we note that for $n$ such that $l_p(n) > \sqrt[4]{x}$, we must have $n^\ast$ prime. Indeed, since
$l_p(n^\ast) > l_p(n)$, for $n^\ast$ not prime we immediately obtain
\[n^\ast > l_p(n)^2 > \sqrt{x} > Y,\]
which contradicts the choice of $n^\ast$. On the other hand, for any $n^\ast$ prime we have
\[n = n^\ast l_p(n) > l_p(n)^2 > Y,\]
giving us a special leaf. \\
Thus the leaves of the first two kinds are of the form $n=pq > Y$, with ${\sqrt[4]{x} < p < q \leq Y}$,
with their contribution being
\begin{align*}
 &\sum_{\sqrt[4]{x} < p \leq Y}\sum_{p < q \leq Y}\mu(pq)f(pq)\phi\left(\frac{x}{pq},\pi(p)-1\right)  \\
=&\sum_{\sqrt[4]{x} < p \leq Y}f(p)\sum_{p < q \leq Y}f(q)\phi\left(\frac{x}{pq},\pi(p)-1\right).
\end{align*}

\subsection{Special Leaves \romannum{1}}

For $\sqrt[3]{x} < p < q \leq Y$ we immediately notice that
\[pq > \sqrt[3]{x^2} \text{, and } \frac{x}{pq} < \sqrt[3]{x} < p,\]
thus $\phi\left(\frac{x}{pq},\pi(p)-1\right) = 1$. Their contribution is then
\begin{align*}S_{c1}
    &= \sum_{\sqrt[3]{x} < p \leq Y}f(p)\sum_{p < q \leq Y}f(q)                         \\
    &= \sum_{\sqrt[3]{x} < p \leq Y}f(p)\left(F(Y) - F(p)\right)                        \\
    &= F(Y)\sum_{\sqrt[3]{x} < p \leq Y}f(p) - \sum_{\sqrt[3]{x} < p \leq Y}f(p)F(p)    \\
    &= F(Y)\left(F(Y)-F(\sqrt[3]{x})\right) - \sum_{\sqrt[3]{x} < p \leq Y}f(p)F(p),    \\
\end{align*}
and we can compute this immediately after sieving $[1,Y]$

\subsection{Special Leaves \romannum{2}.1}

Now we consider the leaves with $\sqrt[4]{x} < p \leq \sqrt[3]{x}$, and $q > \nicefrac{x}{p^2}$.
This gives us $p^2 > \nicefrac{x}{q} \geq \nicefrac{x}{Y}$, and $\frac{x}{pq} < p$. Thus we have again
$\phi\left(\frac{x}{pq},\pi(p)-1\right) = 1$, so their contribution is
\begin{align*}S_{c21}
    &= \sum_{\sqrt{\nicefrac{x}{Y}} < p \leq \sqrt[3]{x}}f(p)\sum_{\nicefrac{x}{p^2} < q \leq Y}f(q)            \\
    &= \sum_{\sqrt{\nicefrac{x}{Y}} < p \leq \sqrt[3]{x}}f(p)\left(F(Y) - F\left(\frac{x}{p^2}\right)\right)    \\
    &= F(Y)\left(F(\sqrt[3]{x}) - F\left(\sqrt{\frac{x}{Y}}\right)\right)
            -\sum_{\sqrt{\nicefrac{x}{Y}} < p \leq \sqrt[3]{x}}f(p)F\left(\frac{x}{p^2}\right).
\end{align*} \\
Since $\nicefrac{x}{p^2} \leq Y$, we may assume that $F\left(\nicefrac{x}{p^2}\right)$ was precomputed. \\

\subsection{Special Leaves \romannum{2}.2}

We consider the leaves with $\sqrt[4]{x} < p \leq \sqrt[3]{x}$, and $q \leq \nicefrac{x}{p^2}$.
For $\phi\left(\frac{x}{pq},\pi(p)-1\right)$ we want the terms not divisible by primes below $p$,
and we have $p < \frac{x}{pq} < \sqrt{x} < p^2$. Thus the terms we want are exactly $1$ and the prime numbers
in the interval $\left[p, \frac{x}{pq}\right]$, giving us
\begin{align*}
\phi\left(\frac{x}{pq}, \pi(p)-1\right)
    &= 1 + \sum_{\substack{p \leq r \leq\frac{x}{pq} \\ r \text{ prime }}}f(r)    \\
    &= 1 + F\left(\frac{x}{pq}\right) - F(p-1).
\end{align*}
Thus the total contribution of these nodes is
\begin{align*}
S_{c22}
    &= \sum_{\sqrt[4]{x} < p \leq \sqrt[3]{x}}f(p)
        \sum_{p < q \leq \min{(\nicefrac{x}{p^2}, Y)}}f(q)\left(1 + F\left(\frac{x}{pq}\right) - F(p-1)\right)  \\
    &= \sum_{\sqrt[4]{x} < p \leq \sqrt[3]{x}}f(p)\left(1 - F(p-1)\right)
        \sum_{p < q \leq \min{(\nicefrac{x}{p^2}, Y)}}f(q)                                                      \\
    &+ \sum_{\sqrt[4]{x} < p \leq \sqrt[3]{x}}f(p)
        \sum_{p < q \leq \min{(\nicefrac{x}{p^2}, Y)}}f(q)F\left(\frac{x}{pq}\right).
\end{align*}
We note, that $Y \leq \nicefrac{x}{p^2}$ when $p^2 \leq \nicefrac{x}{Y}$. We start by using this to split
the first sum above into the sums without conditions:
\begin{align*}
    &\sum_{\sqrt[4]{x} < p \leq \sqrt{\nicefrac{x}{Y}}}f(p)\left(1 - F(p-1)\right)
        \sum_{p < q \leq Y}f(q).                                                                            \\
    &= \sum_{\sqrt[4]{x} < p \leq \sqrt{\nicefrac{x}{Y}}}f(p)\left(1 - F(p-1)\right)\left(F(Y)-F(p)\right). \\
    &\sum_{\sqrt{\nicefrac{x}{Y}} < p \leq \sqrt[3]{x}}f(p)\left(1 - F(p-1)\right)
        \sum_{p < q \leq \nicefrac{x}{p^2}}f(q)                                                             \\
    &= \sum_{\sqrt{\nicefrac{x}{Y}} < p \leq \sqrt[3]{x}}f(p)\left(1 - F(p-1)\right)
        \left(F\left(\frac{x}{p^2}\right)-F(p)\right).
\end{align*}
That gives us
\[\sum_{\sqrt[4]{x} < p \leq \sqrt[3]{x}}f(p)\left(1 - F(p-1)\right)
    \sum_{p < q \leq \min{(\nicefrac{x}{p^2}, Y)}}f(q) = V_1 + V_2 - V_3,\]
with
\begin{align*}
V_1 &= \sum_{\sqrt[4]{x} < p \leq \sqrt{\nicefrac{x}{Y}}}f(p)\left(1 - F(p-1)\right)F(Y)                \\
    &= F(Y)\sum_{\sqrt[4]{x} < p \leq \sqrt{\nicefrac{x}{Y}}}f(p)\left(1 - F(p-1)\right)                \\
    &= F(Y)\left(F\left(\sqrt{\frac{x}{Y}}\right) - F(\sqrt[4]{x})\right) -
            F(Y)\sum_{\sqrt[4]{x} < p \leq \sqrt{\nicefrac{x}{Y}}}f(p)F(p-1)                            \\
V_2 &= \sum_{\sqrt{\nicefrac{x}{Y}} < p \leq \sqrt[3]{x}}f(p)
            \left(1 - F(p-1)\right)F\left(\frac{x}{p^2}\right)                                          \\
    &= \sum_{\sqrt{\nicefrac{x}{Y}} < p \leq \sqrt[3]{x}}f(p)F\left(\frac{x}{p^2}\right) -
            \sum_{\sqrt{\nicefrac{x}{Y}} < p \leq \sqrt[3]{x}}f(p)F(p-1)F\left(\frac{x}{p^2}\right)     \\
V_3 &= \sum_{\sqrt[4]{x} < p \leq \sqrt[3]{x}}f(p)\left(1 - F(p-1)\right)F(p)                           \\
    &= \sum_{\sqrt[4]{x} < p \leq \sqrt[3]{x}}f(p)F(p) - \sum_{\sqrt[4]{x} < p \leq \sqrt[3]{x}}f(p)F(p-1)F(p)
\end{align*}
As before we note that all the values of $F$ has the argument in $[1,Y]$.

We proceed similarly for the second sum:
\[\sum_{\sqrt[4]{x} < p \leq \sqrt[3]{x}}f(p)
    \sum_{p < q \leq \min{(\nicefrac{x}{p^2}, Y)}}f(q)F\left(\frac{x}{pq}\right) = W_1 + W_2,\]
where
\begin{align*}
W_1 &= \sum_{\sqrt[4]{x} < p \leq \sqrt{\nicefrac{x}{Y}}}f(p)
            \sum_{p < q \leq Y}f(q)F\left(\frac{x}{pq}\right)                   \\
W_2 &= \sum_{\sqrt{\nicefrac{x}{Y}} < p \leq \sqrt[3]{x}}f(p)
            \sum_{p < q \leq \nicefrac{x}{p^2}}f(q)F\left(\frac{x}{pq}\right)   \\
\end{align*}
We note that we can proceed in a similar way to the calculation of $P_2$: as we sieve the block, for primes
$p \in \left(\sqrt[4]{x}, \sqrt[3]{x}\right]$ we find bounds on $q$, so that $F\left(\frac{x}{pq}\right)$ is in the
current block. Noting that $q \leq Y$, unlike $P_2$, we do not need to sieve the resulting interval. Further,
$\frac{x}{pq} \leq \sqrt{x}$, give us the upper bound on blocks we need to consider. \\

\subsection{Special Leaves: Bringing it all together}

As we have seen above, all the terms in $S_{c1}$, $S_{c21}$, and $V_1,V_2,V_3$ can be precomputed.
We should try to combine these sums. First of all, the constant term is
\begin{align*}
&F(Y)\left(F(Y) - F(\sqrt[3]{x})\right) +
    F(Y)\left(F(\sqrt[3]{x}) - F\left(\sqrt{\frac{x}{Y}}\right)\right) +
    F(Y)\left(F\left(\sqrt{\frac{x}{Y}}\right) - F(\sqrt[4]{x})\right)      \\
&= F(Y)\left(F(Y) - F(\sqrt[4]{x})\right).
\end{align*}
Then we combine the terms involving $f(p)F(p)$ to get: $-\sum_{\sqrt[4]{x} < p \leq Y}f(p)F(p)$. \\
Further, the terms involving $f(p)F\left(\nicefrac{x}{p^2}\right)$ are annihilated. \\
Thus
\begin{align*}S_1
    &= S_{c1}+S_{c21}+V_1+V_2+V_3                                                                   \\
    &= F(Y)\left(F(Y) - F(\sqrt[4]{x})\right) - \sum_{\sqrt[4]{x} < p \leq Y}f(p)F(p)               \\
    &- F(Y)\sum_{\sqrt[4]{x} < p \leq \sqrt{\nicefrac{x}{Y}}}f(p)F(p-1)                             \\
    &- \sum_{\sqrt{\nicefrac{x}{Y}} < p \leq \sqrt[3]{x}}f(p)F(p-1)F\left(\frac{x}{p^2}\right)      \\
    &+ \sum_{\sqrt[4]{x} < p \leq \sqrt[3]{x}}f(p)F(p-1)F(p),
\end{align*}
and, as noted above, all the values used can be computed once $[1,Y]$ is sieved. \\
We set $S_2 = W_1 + W_2$, which should be updated per-block in $[1,\sqrt{x}]$.
And we set $S_3$ to be the contribution of the special leaves of the third class,
which should use the algorithm described above to update the values of $\phi$. \\

Thus, to summarize:
\begin{samepage}\begin{enumerate}
    \item We calculate the contribution of ordinary leaves.
    \item We calculate $S_1$.
    \item For each block in $[1,\sqrt{x}]$ we update $S_2$.
    \item For each block we update the table of $\phi(n,b)$ for $b \leq \sqrt[4]{x}$ and use it to calculate $S_3$.
\end{enumerate}\end{samepage}

For $S_1$ we note that, to allocate less memory, all the terms apart from the sum of
$f(p)F(p-1)F\left(\nicefrac{x}{p^2}\right)$ can be calculated per block. And for the latter we need to only store
$F(p-1)$ for primes in $\left(\sqrt{\nicefrac{x}{Y}}, \sqrt[3]{x}\right]$, noting that this will be
available before we calculate $F\left(\nicefrac{x}{p^2}\right)$, so we can, again, calculate this sum per-block. \\
We can lessen the memory usage further. Note that $p \leq \nicefrac{x}{p^2}$. Suppose we have calculated $F(p-1)$ for
the largest $p$ not greater than $\sqrt[3]{x}$. Afterwards we will encounter the values $\nicefrac{x}{p^2}$ in the
order of $p$ descending; thus after processing $p_k$ we subtract $f(p_{k-1})$ from the accumulator and continue.

\ifstandalone
\bibliography{master}
\bibliographystyle{amsplain}
\fi

\section{Parallelizing Computations}

We note that a lot of computations can be easily parallelized. We start with $P_2$ computation.

\subsection{Computation of $P_2$}

As before, we proceed in blocks $I=[l,l+B)$. Let's write $F_I$ for ``in block'' summatory function of $f$, \ie
$F_I(x) = \sum_{l \leq p \leq x}f(p)$, assuming $x \in I$. Then the contribution of $I$ into $S_2$ would be
\begin{align*}S_{2I}
    &= \sum_{l \leq p < l+B}f(p)F(p-1) = \sum_{l \leq p < l+B}f(p)\left(F(l-1) + F_I(p-1)\right)    \\
    &= F(l-1)\sum_{l \leq p < l+B}f(p) + \sum_{l \leq p < l+B}f(p)F_I(p-1).
\end{align*}
Thus to compute $S_2$ in parallel, we split it into two parts: first of all we sum the values of
$\sum_{l \leq p < l+B}f(p)F_I(p-1)$ per block, and also we compute $F_I(l+B-1)$ per block. The latter values can be
used to reconstruct $F(l-1)$ in order. \\
We proceed similarly for $S_1$
\begin{align*}S_{1I}
    &= \sum_{l \leq \nicefrac{x}{p} < l+B}f(p)F\left(\frac{x}{p}\right) =
        \sum_{l \leq \nicefrac{x}{p} < l+B}f(p)\left(F(l-1) + F_I\left(\frac{x}{p}\right)\right)    \\
    &= F(l-1)\sum_{l \leq \nicefrac{x}{p} < l+B}f(p) +
        \sum_{l \leq \nicefrac{x}{p} < l+B}f(p)F_I\left(\frac{x}{p}\right),
\end{align*}
and once again, the second sum is fully computed per block, and the first one is reconstructed. \\

\subsection{Computation of $\phi$}

For $S_1$ we note that we have sums involving products $F(p-1)F\left(\nicefrac{x}{p^2}\right)$ and $F(p-1)F(p)$ which
is hard to handle in the same manner as above. Thus we should handle it separately, either by going in blocks through
$[1,Y]$, or, as we sieve $[1,Y]$ anyway, we can precompute all the values of $F(n)$ in this interval and then run in
blocks in parallel, as then we won't have any interdependency. \\
For $S_2$ we proceed similar to the calculation of $P_2$:
\begin{align*}S_{2I}
    &= \sum_{\nicefrac{x}{pq} \in I} f(p)f(q)F\left(\frac{x}{pq}\right) =
        \sum_{\nicefrac{x}{pq} \in I} f(p)f(q)\left[F(l-1) + F_I\left(\frac{x}{pq}\right)\right]    \\
    &= F(l-1)\sum_{\nicefrac{x}{pq} \in I} f(p)f(q) +
        \sum_{\nicefrac{x}{pq} \in I} f(p)f(q)F_I\left(\frac{x}{pq}\right).
\end{align*} \\
For $S_3$ we write
\begin{gather*}
G_I(x,k) = \sum_{\substack{n \in I, n \leq x \\ (n,m_k) = 1}}f(n),  \\
g(n,b) = \mu(n)f(n)f(p_{b+1}),
\end{gather*}
and we have
\begin{align*}S_{3I}
    &= -\sum_{\nicefrac{x}{n^\ast p_{b+1}} \in I}g(n^\ast, p_{b+1})
        \left[\phi(l-1,b) + G_I\left(\frac{x}{n^\ast p_{b+1}},b\right)\right]               \\
    &= -\sum_b \phi(l-1,b) \sum_{\nicefrac{x}{n^\ast p_{b+1}} \in I}g(n^\ast, p_{b+1})
        -\sum_{\nicefrac{x}{n^\ast p_{b+1}} \in I}g(n^\ast, p_{b+1})
            G_I\left(\frac{x}{n^\ast p_{b+1}},b\right)
\end{align*}

We note that we touch only ``trivial'' parallelization concerns. We refer to \cite{KimWalischCount} for an in-deep
discussion of the possible optimizations, as here we are mostly interested in maths.

\ifstandalone
\bibliography{master}
\bibliographystyle{amsplain}
\fi

\section{Computational Analysis}

As we do the analysis to find the guidelines to select the good values for $Y$ and $B$, we should ignore the
impact of the multi-precision arithmetic (if it is used), so we assume that multiplications, additions and the
calculation of $f$ is $\mathcal{O}(1)$.

\subsection{Cost of Sieving}

Our main subtask will be to sieve the interval of length $N$ with primes up to $Y$. The complexity can be estimated as
\[N\sum_{p \leq Y}\frac{1}{p} \sim N\log\log{Y}.\]

We start with sieving $[1,Y]$, and then we fully sieve $[Y,\frac{x}{Y}]$ in blocks of size $B$ using only primes
up to $Y$. The total complexity of the sieving step is then
\[\mathcal{O}\left(\frac{x}{Y}\log\log{Y}\right).\]

\subsection{Cost of $P_2$}
For $S_1$ we need to sieve disjoint subintervals of $[\sqrt{x}, \nicefrac{x}{Y})$ to find $q$ for $\nicefrac{x}{q}$
to be in the current block. We may estimate this as the time to fully sieve $[\sqrt{x}, \nicefrac{x}{Y})$ giving us
\[\mathcal{O}\left(\frac{x}{Y}\log\log{Y}\right).\]
Then we sum $\pi(\sqrt{x})-\pi(Y)$ terms giving us an estimate of
\[\mathcal{O}\left(\pi(\sqrt{x})-\pi(Y)\right) = \mathcal{O}\left(\frac{\sqrt{x}}{\log{x}}\right).\]
Since $\sqrt{x} \leq \nicefrac{x}{Y}$ we have the cost of $P_2$ to be
\[\mathcal{O}\left(\frac{x}{Y}\log\log{Y}\right).\]

\subsection{Cost of $\phi$: $S_1$}

In here we have simple sums, so the total complexity depends just on the number of terms, and can be approximated by
\[\mathcal{O}\left(\pi(Y)\right) = \mathcal{O}\left(\frac{Y}{\log{Y}}\right).\]

\subsection{Cost of $\phi$: $S_2$}

To estimate the complexity of computing $W_1$ we need to count the number of terms in the sum. This gives us
\begin{align*}
&   \left(\pi\left(Y\right) - \pi\left(\sqrt[4]{x}\right)\right) + \dots +
        \left(\pi\left(Y\right) - \pi\left(\sqrt{\frac{x}{Y}}\right)\right)                             \\
&=  \frac{\pi\left(\sqrt{\frac{x}{Y}}\right) - \pi\left(\sqrt[4]{x}\right)}{2}
        \left[2\pi\left(Y\right) - \left(\pi\left(\sqrt{\frac{x}{Y}}\right) + \pi\left(\sqrt[4]{x}\right)\right)\right]
\end{align*}
First of all we note that $\sqrt[3]{x} \leq Y < \sqrt{x}$ hence we have
\[\sqrt[4]{x} < \sqrt{\frac{x}{Y}} \leq \sqrt[3]{x}.\]
Thus the sum above can be approximated by
\[\mathcal{O}\left(\frac{Y\sqrt[3]{x}}{\log^2{x}}\right).\]

Similarly for $W_2$ we have $\sqrt[3]{x} \leq \nicefrac{x}{p^2} < Y$ and
\[\mathcal{O}\left(\sum_{\sqrt{\nicefrac{x}{Y}} < p \leq \sqrt[3]{x}} \pi\left(\frac{x}{p^2}\right)\right)
    = \mathcal{O}\left(\pi\left(\sqrt[3]{x}\right)\right) \mathcal{O}\left(\pi\left(Y\right)\right)
    = \mathcal{O}\left(\frac{Y\sqrt[3]{x}}{\log^2{x}}\right).
\]

\subsection{Cost of $\phi$: $S_3$}
The process of sieving blocks is made complicated here by the Fenwick Tree. For each node of interest we need
to calculate the prefix sum, giving us $\mathcal{O}(\log{B})$ complexity.
To estimate the number of the nodes we need to process we note that they have form $(n^\ast p_{b+1}, b)$
with $n^\ast \leq Y$ and $b \leq \pi\left(\sqrt[4]{x}\right)$, giving us a total of
\[\mathcal{O}\left(Y\pi\left(\sqrt[4]{x}\right)\log{B}\right)\]
Also, we need to account for updating the Fenwick tree. We note that we touch every value at most once: first time
when we initialize the tree for a block, and a second time when we ``strike out'' this value\footnotemark,
thus giving us an estimation on the extra work to be
\[\mathcal{O}\left(\frac{x}{Y}\log{B}\right).\]
\footnotetext{We should note that this is a very crude estimation as we strike out the multiples of the first
$K$ primes without any extra processing.}

\subsection{Total Cost}

We see that the total cost is about
\[\mathcal{O}\left(\frac{x}{Y}\log\log{Y} + \frac{Y}{\log{Y}}
    + \frac{Y\sqrt[3]{x}}{\log^2{x}}
    + \frac{x}{Y}\log{B} + \frac{Y\sqrt[4]{x}\log{B}}{\log{x}}
\right).\]
Ignoring $B$ for now (we can approximate it with $\log{B} \sim \log{Y}$) we note that setting
$Y=\mathcal{O}\left(\sqrt[3]{x}\log^3{x}\right)$ will give us the total estimate of
\[\mathcal{O}\left(\frac{x^{\frac{2}{3}}}{\log^2{x}}\right).\]

Following \cite{KimWalischSum} we should set $Y = \alpha\sqrt[3]{x}$ with
\[\alpha = a\log^3{x}+b\log^2{x}+c\log{x}+d,\]
and find the optimal values for $a,b,c,d$ empirically.

\subsection{Finding $\alpha$}

As noted in \cite{Silva},
``changes of $\pm25\%$ around the optimal value of $\alpha$ did not increase the execution time by more than $3\%$'',
which agrees with our experiments. Thus for a given $x$ we fix the interval $[\alpha_0,\alpha_1]$ and calculate
$F(x)$ using a different $\alpha$, finding the one that gives us the fastest running time. We repeat this for several
values of $x$ and then we find the best fitting values. We must say that this process is pretty ``noisy''
thus we did several runs of it. Even then, from our experiments, the optimal value of $\alpha$ does not change
considerably with changing the power exponent, thus we can have one formula to calculate $\alpha$ for all.
The results we've got on the author's computer are as follows:

\begin{center}
\begin{tabular}{|l|c|c|c|c|c|c|c|c|c|c|}
\hline $x$ & $10^9$ & $10^{10}$ & $5 \times 10^{10}$ & $10^{11}$ & $5 \times 10^{11}$ & $10^{12}$ &
    $5 \times 10^{12}$ & $10^{13}$ & $5 \times 10^{13}$ & $10^{14}$ \\
\hline $\alpha$ & 3 & 3.5 & 5 & 6 & 7 & 8 &
    8.5 & 10 & 12 & 13                                              \\
\hline
\end{tabular}
\end{center}

Fitting the values to the data we obtain the formula
\[\alpha \approx 0.000681 \log{x}^3 - 0.011846  \log{x}^2 + 0.044074 \log{x} + 0.988365.\]

\ifstandalone
\bibliography{master}
\bibliographystyle{amsplain}
\fi

\section{Numerical Results}

\noindent We give some values of the function $F_n(x) = \sum_{p \leq x}p^n$. \\
We have checked the correctness of our algoritms in several ways. First of all, for $n=0$ and $n=1$ we have compared
the results to the known values: for $n=0$ Wikipedia has an extensive table and for $n=1$ we have used the values
from \cite{KimWalischSum}. For ``smaller'' values of $x$ we have checked the results against Pari/GP. We have also
computed $F_n(x)$ and $F_n(x+\varepsilon)$, and checked that these agree by sieving $[x,x+\varepsilon]$ and calculating
the value of $F_n(x+\varepsilon)$ from $F_n(x)$ and now known primes. We have also tested extensively for the small $x$,
using the variation of parameters: $Y$, $B$, and $K$. \\

We note that we have used 256-bit arithmetics for our calculations, so the table for $F_4(x)$ is shorter due to overflows
with larger values of $x$.

\vspace{7mm}
\begin{table}[htp!]
\caption{Values of $F_2(x)$}
\centering
\begin{tabular}{|r|r|} \hline
$10^{10}$ & 14692485666215945973239505690 \\ \hline
$5 \times 10^{10}$ & 1714863031171407826702942323341 \\ \hline
$10^{11}$ & 13338380640732671147186590712800 \\ \hline
$5 \times 10^{11}$ & 1566398144419578032981266419280441 \\ \hline
$10^{12}$ & 12212907966177661747436156685876997 \\ \hline
$5 \times 10^{12}$ & 1441593988892141564900337100187358316 \\ \hline
$10^{13}$ & 11262617785640702236670513970349205634 \\ \hline
$5 \times 10^{13}$ & 1335210125295770298473184342618020082018 \\ \hline
$10^{14}$ & 10449549945144268110573967892555485354493 \\ \hline
$5 \times 10^{14}$ & 1243450253668811479272045017247069947749359 \\ \hline
$10^{15}$ & 9745981795365753183493378490092915742101696 \\ \hline
$5 \times 10^{15}$ & 1163492503926172589836028128116501190925639911 \\ \hline
$10^{16}$ & 9131187419861160902346450308274850949251333488 \\ \hline
$5 \times 10^{16}$ & 1093197297594496716923810873727618016729821706326 \\ \hline
$10^{17}$ & 8589360822439890567209125328673103991944617500636 \\ \hline
\end{tabular}
\end{table}

\vspace{7mm}
\begin{table}[htp!]
\caption{Values of $F_3(x)$}
\centering
\begin{tabular}{|r|r|} \hline
$10^{10}$ & 109780001885333601058528339379120755908 \\ \hline
$5 \times 10^{10}$ & 64082046820723451487075613889900198582245 \\ \hline
$10^{11}$ & 996973732171667396998099396013424430102364 \\ \hline
$5 \times 10^{11}$ & 585523008895664909752699795716553990313510861 \\ \hline
$10^{12}$ & 9131183180200496139738672227721825939508051079 \\ \hline
$5 \times 10^{12}$ & 5390112237650538623349869931020562476048206694122 \\ \hline
$10^{13}$ & 84227641426129569665463994634958535318856709486404 \\ \hline
$5 \times 10^{13}$ & 49934458175987733033055840585837321485798098600691128 \\ \hline
$10^{14}$ & 781635768313974029776529273260033374787234114404747147 \\ \hline
$5 \times 10^{14}$ & 465116388359752442960756526955719292849606973879669769221 \\ \hline
$10^{15}$ & 7291408725599572782932502326635718783267041072449283373284 \\ \hline
$5 \times 10^{15}$ & 4352798262365871370931577893185982626698534306936825507930393 \\ \hline
$10^{16}$ & 68325370066732554478047827274497012925608706449875092937488310 \\ \hline
$5 \times 10^{16}$ & 40904060647664055185674494652305281362155344363027361735219374850 \\ \hline
$10^{17}$ & 642800454307687984344682304535086826246739426032801572692032867900 \\ \hline
\end{tabular}
\end{table}

\vspace{7mm}
\begin{table}[htp!]
\caption{Values of $F_4(x)$}
\centering
\begin{tabular}{|r|r|} \hline
$10^{10}$ & 876279913324387539183894015044723229219045342750 \\ \hline
$5 \times 10^{10}$ & 2557935537494958740759502417763799692339226020428977 \\ \hline
$10^{11}$ & 79596284512301003834661995024051166148172171338471852 \\ \hline
$5 \times 10^{11}$ & 233763531449591824849643583298320925979435007893090787181 \\ \hline
$10^{12}$ & 7291405674369073069761776122154275329527216556249514408897 \\ \hline
$5 \times 10^{12}$ & 21522735063683505999117268748154742855602976344249380004455248 \\ \hline
$10^{13}$ & 672670357329861606491888367274980639337624917176989202367080838 \\ \hline
$5 \times 10^{13}$ & 1994144875417964028004350412269602745929105301928382366911602202806 \\ \hline
$10^{14}$ & 62431856891826059995323082674283840302039995668819924057075933466801 \\ \hline
$5 \times 10^{14}$ & 185766273878955147609467369321682910437837624744901019253719657653093899 \\ \hline
$10^{15}$ & 5824519129976593880511325120158491808935524746219691118247177308917094908 \\ \hline
\end{tabular}
\end{table}

\ifstandalone
\bibliography{master}
\bibliographystyle{amsplain}
\fi

\FloatBarrier
\bibliographystyle{amsplain}
\bibliography{master}

\providecommand{\bysame}{\leavevmode\hbox to3em{\hrulefill}\thinspace}
\providecommand{\MR}{\relax\ifhmode\unskip\space\fi MR }
\providecommand{\MRhref}[2]{%
  \href{http://www.ams.org/mathscinet-getitem?mr=#1}{#2}
}
\providecommand{\href}[2]{#2}
\begin{thebibliography}{1}

\bibitem{DelegliseRivat}
M.~Deleglise and J.~Rivat, \emph{Computing $\pi(x)$: the {Meissel}, {Lehmer},
  {Lagarias}, {Miller}, {Odlyzko} method}, Mathematics of Computation
  \textbf{65} (1996), no.~213, 235--246 (en).

\bibitem{LagariasMillerOdlyzko}
J.~C. Lagarias, V.~S. Miller, and A.~M. Odlyzko, \emph{Computing $\pi(x)$: the
  {Meissel}-{Lehmer} method}, Mathematics of Computation \textbf{44} (1985),
  no.~170, 537--537.

\bibitem{LagariasOdlyzkoAnal}
J.C Lagarias and A.M Odlyzko, \emph{Computing $\pi(x)$: An analytic method},
  Journal of Algorithms \textbf{8} (1987), no.~2, 173--191.

\bibitem{Lehmer}
D.~H. Lehmer, \emph{On the exact number of primes less than a given limit},
  Illinois Journal of Mathematics \textbf{3} (1959), no.~3, 381--388 (en).

\bibitem{Silva}
T.~E.~O. Silva, \emph{Computing $\pi$(x): the combinatorial method}, Revista do
  DETUA \textbf{4} (2006), no.~6, 759--768.

\bibitem{KimWalischCount}
Kim Walisch, \emph{https://github.com/kimwalisch/primecount}, 2021.

\bibitem{KimWalischSum}
\bysame, \emph{https://github.com/kimwalisch/primesum}, 2021.

\end{thebibliography}

\end{document}